\documentclass[twocolumn]{autart} 
\usepackage[numbers,sort&compress]{natbib}  
\usepackage{color}
\usepackage{amsmath,amssymb,amsfonts}
\usepackage{xpatch}
\xpatchcmd{\proof}{\itshape}{\prooflabelfont}{}{}
\newcommand{\prooflabelfont}{\bfseries}
\usepackage{subcaption}
\usepackage{empheq}
\usepackage{stfloats}
\makeatletter \let\cl@part\relax \makeatother
\usepackage{algorithm}
\usepackage{algpseudocode}
\usepackage{graphicx}
\usepackage{booktabs}
\usepackage{arydshln}
\usepackage{multirow}
\usepackage{wrapfig}
\usepackage[capitalise]{cleveref}

\newtheorem{myrem}{Remark}
\newtheorem{mylem}{Lemma}
\newtheorem{mythr}{Theorem}

\newtheorem{myass}{Assumption}

\newtheorem{mypro}{Proposition}

\begin{document}

\begin{frontmatter}

\title{Data-Driven Structured Controller Design Using the  Matrix S-Procedure
% \thanksref{footnoteinfo}
} % Title, preferably not more 
                                                % than 10 words.

% \thanks[footnoteinfo]{}

\author[hkust]{Zhaohua Yang}\ead{zyangcr@connect.ust.hk},
~~    
\author[hkust]{Yuxing Zhong}\ead{yuxing.zhong@connect.ust.hk},
~~         
\author[alberta]{Nachuan Yang}\ead{nachuan1@ualberta.ca},
~~ 
\author[hkust]{Xiaoxu Lyu}\ead{eelyuxiaoxu@ust.hk},
~~  
\author[hkust]{Ling Shi}\ead{eesling@ust.hk}

\address[hkust]{Department of Electronic and Computer Engineering, The Hong Kong University of Science and Technology,\\ 
Clear Water Bay, Hong Kong, China}
\address[alberta]{Department of Electrical and Computer Engineering, University of Alberta,\\
Edmonton, Alberta T6G 1H9, Canada}
\maketitle

\begin{abstract}
This paper focuses on the data-driven optimal structured controller design for discrete-time linear time-invariant (LTI) systems, considering both the $H_2$ performance and the $H_\infty$ performance. Specifically, we consider three scenarios: (i) the model-based structured control, (ii) the data-driven unstructured control, and (iii) the data-driven structured control. For the $H_2$ performance, we primarily investigate cases (ii) and (iii), since case (i) has been extensively studied in the literature. For the $H_\infty$ performance, all three scenarios are considered. For the structured control, we introduce a linearization technique that transforms the original nonconvex problem into a semidefinite programming (SDP) problem. Based on this transformation, we develop an iterative linear matrix inequality (ILMI) algorithm. For the data-driven control, we describe the set of all possible system matrices that can generate the sequence of collected data. Additionally, we propose a sufficient condition to handle all possible system matrices using the matrix S-procedure. The data-driven structured control is followed by combining the previous two cases. We compare our methods with those in the existing literature and demonstrate our advantages via several numerical simulations.
\end{abstract}

\end{frontmatter}

\section{Introduction}
Networked control systems (NCS) have found numerous applications in real-world scenarios. In such systems, the feedback controller establishes a communication channel connecting the system's outputs and inputs. Consequently, introducing sparsity into the controller is essential, as a centralized control strategy can impose a substantial communication burden in large-scale systems. This motivates us to design a controller with a specified sparse structure.

Designing a controller with structural constraints is challenging due to its NP-hardness \cite{blondel1997np}. The past decades have witnessed various literature working on this problem \cite{polyak2013lmi,lin2011augmented,jovanovic2016controller,fardad2014design,yang2024log}. The methods can be primarily classified into LMI-based and gradient-based approaches. Polyak et al. \cite{polyak2013lmi} showed that designing a row-sparse controller is not conservative, as post-multiplication preserves the row structure. In addition, Jovanović et al. \cite{jovanovic2016controller} designed a structured controller by restricting the Lyapunov matrix to be diagonal. However, this method can be overly conservative and may lead to infeasibility. Fardad et al. \cite{fardad2014design}, on the other hand, developed ILMI algorithms to compute the structured controller. In the seminal work of the gradient-based method, Lin et al. \cite{lin2011augmented} leveraged the augmented Lagrangian method to design a structured controller to optimize the system's $H_2$ performance.  Recently, Yang et al. \cite{yang2024log} designed a new algorithm based on the interior point method. 
% However, most of the existing literature assumes the system model is known \textit{a priori}, with very few studies addressing the model-free case, which motivates our research.
The aforementioned literature assumes the system model is known \textit{a priori}. However, in reality, the dynamics of many systems are often inaccessible, which motivates research on data-driven control.

The past decades have witnessed extensive research on data-driven control, which operates under the assumption that the system is unknown and only a sequence of sampled data is accessible. A fundamental tool for analysis in this area is \textit{Willems et al.'s fundamental lemma \cite{willems2005note}}, which provides criteria for determining whether the collected data is sufficiently informative to identify system dynamics. Various works have emerged based on this lemma. In the seminal work \cite{de2019formulas}, the authors examined robust stabilization with noise satisfying a signal-to-noise ratio constraint. Berberich et al. \cite{berberich2020robust} extended these results and considered the $H_\infty$ performance. Van Waarde et al. \cite{van2020data} further conducted a thorough investigation into the informativity of data for various control problems. Other tools, such as the matrix S-lemma and Petersen's lemma, are also employed to develop direct approaches for the robust case. Van Waarde et al. \cite{van2020noisy} first expressed the uncertainty set of system matrices consistent with data using the quadratic matrix inequality (QMI), and studied robust stabilization as well as $H_2$ and $H_\infty$ control using the matrix S-lemma. Furthermore, Bisoffi et al. \cite{bisoffi2022data,bisoffi2021trade} explicitly characterized this uncertainty set using a matrix ellipsoid and investigated robust stabilization based on Petersen's lemma. They were the first to address instantaneous-bounded noise using the (lossy) matrix S-procedure \cite{bisoffi2021trade}. For sparse control, Eising et al. \cite{eising2022informativity} considered robust stabilization with maximal sparsity, while Miller et al. \cite{miller2024data} designed a structured controller aimed at optimizing the $H_2$ performance. {\color{blue}However, \cite{eising2022informativity} lacked discussions on system performance, and the methods in \cite{miller2024data} designed a structured controller only within a convex subset of the nonconvex feasibility region. These limitations motivate us to propose more advanced approaches for exploring the nonconvex feasibility region in the design of a structured controller.}

This paper focuses on data-driven optimal structured controller design, considering both the $H_2$ performance and $H_\infty$ performance. The main contributions are summarized as follows
\begin{enumerate}
    \item {\color{blue}We propose a unified data-driven structured controller design framework:\\
    (i) Unlike existing literature \cite{van2020noisy,miller2024data} that considered energy-bounded noise, we focus on instantaneous-bounded noise. We characterize the set of all possible system matrices as the intersection of the sets for each data point. Additionally, we utilize the matrix S-procedure to handle all possible system matrices. This approach avoids energy-based overapproximation for instantaneous-bounded noise \cite{miller2024data}, resulting in reduced conservatism at the cost of increased computation. Furthermore, we establish a novel property of our framework that the performance is monotonically nonincreasing w.r.t. the data length.} \\
  (ii) We develop a novel linearization technique that turns the original nonconvex problem into an SDP problem, which allows expressing the structural constraint explicitly. This motivates an ILMI algorithm to compute a structural controller. we discuss that our method can explore a larger feasibility region than {\color{blue}existing techniques}, possibly leading to a structural controller with higher performance.
  \item In the $H_2$ control case, given the extensive literature on model-based structured control, our primary focus is on data-driven unstructured control and data-driven structured control. Our simulations demonstrate that our methods {\color{blue}consistently achieve improved performance compared with \cite{miller2024data} under the same experiment setting}, aligning with our theoretical results.
  \item In the $H_\infty$ control case, we, for the first time, conduct comprehensive research on the model-based structured control, the data-driven unstructured control, and the data-driven structured control. The effectiveness of our methods is validated through a numerical example.
\end{enumerate}
This paper is organized as follows. In \cref{Preliminaries}, some preliminaries on the $H_2$ performance, $H_\infty$ performance, and structured controllers are provided. In \cref{Data collection and problem formulation}, the data collection mechanism is presented and the problem is formulated. In \cref{H2 control}, three subproblems, model-based structured control, data-driven unstructured control, and data-driven structured control, are discussed. {\color{blue}In \cref{comparison with miller}, we theoretically analyze several fundamental properties of our framework.} \cref{Hinf control} is the $H_\infty$ version of \cref{H2 control}. In \cref{Simulations}, some numerical simulations are conducted to validate the effectiveness of our methods. The paper is concluded in \cref{Conclusion}.

\textit{Notations}:
\ Let $\mathbb{R}^{m\times n}$ and $\mathbb{R}^n$ denote the sets of real matrices of size $m\times n$ and real vectors of size $n$, respectively. For a scalar $\alpha$, $\{\alpha_i\}$ denotes a sequence of scalar. For a matrix $P$, $P^\top$ denotes its transpose, $\operatorname{Tr}(P)$ denotes its trace, and $P_{ij}([P]_{ij})$ denotes its element at $i$-th row
and $j$-th column. The symbol $P\succ 0$ ($P\succeq 0$) means $P$ is positive definite (positive semidefinite). Let $P_{i:j}$ represent the submatrix from $P$'s $i$-th column to $j$-th column. The identity matrix, zero matrices, and all-one matrices are denoted as $I$, $0$, and $\textbf{1}$, with dimensions labeled when necessary. Let $\mathbb{S}^n$ denote the set of symmetric matrices of size $n\times n$. Let $\mathbb{S}^n_+$ ($\mathbb{S}^n_{++}$) denote the set of positive semidefinite matrices (positive definite matrices). Let $\forall$ denote ``for all''. Let $\mathrm{blkdiag}(A_1,A_2,\ldots,A_n)$ denote a block diagonal matrix with $A_i, i=1,2,\ldots,n$ as its blocks. Let $\circ$ represent the elementwise product and $\triangleq$ denote the equality for definition.

\section{Preliminaries}\label{Preliminaries}
In this paper, we consider a discrete-time LTI system
\begin{equation}\label{basic system}
\begin{split}
    x_{k+1} &= Ax_k + Bu_k + Gd_k,\\
    y_k &= Cx_k + Du_k + Hd_k,
\end{split}
\end{equation}
where $x_k\in \mathbb{R}^{n_x}$ is the system state, $u_k\in \mathbb{R}^{n_u}$ is the control input, $y_k\in \mathbb{R}^{n_y}$ is the controlled output, and $d_k\in \mathbb{R}^{n_d}$ is the exogenous disturbance. It is standard to assume that $(A,B)$ is stabilizable. We regulate the system \eqref{basic system} using a state-feedback controller
\begin{equation}
    u_k = Kx_k,
\end{equation}
where $K\in \mathbb{R}^{n_u\times n_x}$ is the control gain to be determined. With this controller, the system \eqref{basic system} can be rewritten as a closed-loop system
\begin{equation}\label{controlled system}
    \begin{split}
        x_{k+1}& = (A+BK)x_k + Gd_k,\\
        y_k &= (C+DK)x_k + Hd_k.
    \end{split}
\end{equation}
The transfer function in the $z$-domain from $d$ to $y$ is given by
\begin{equation}
    T_{yd}(z) = (C+DK)(zI-(A+BK))^{-1}G+H.
\end{equation}
Given the transfer function $T_{yd}(z)$, we denote the $H_2$ norm as $||T_{yd}(z)||_2$ and the $H_\infty$ norm as $||T_{yd}(z)||_\infty$. It is well known that finding a suboptimal $H_2$ or $H_\infty$ controller can be achieved by solving a feasibility LMI problem. To facilitate later discussions, we provide relevant results below.

\begin{mylem}[\cite{zhou1996robust}]
    Assume $H=0$. The $H_2$ norm of the system \eqref{controlled system} satisfies $||T_{yd}(z)||_2\le\gamma$ if and only if there exists matrices $P\in\mathbb{S}^{n_x}_{++}, Q\in\mathbb{S}^{n_y}_{++}$ such that
    \begin{subequations}\label{H2 feasiblity}
    \begin{align}
    &\begin{bmatrix}
	P		&(A+BK)P		&G		\\
	P(A+BK)^\top		&P		&0		\\
	G^\top		&0		&I		\\
    \end{bmatrix} \succ 0, \label{H2 only involved AB}\\
    & \left[ \begin{matrix}\centering
	Q		&(C+DK)P		\\
	P(C+DK)^\top		&P\\
    \end{matrix} \right] \succ 0, \operatorname{Tr}(Q)\le \gamma^2.
    \end{align}
    \end{subequations}
\end{mylem}

\begin{mylem}[\cite{gahinet1994linear}]
    The $H_\infty$ norm of the system \eqref{controlled system} satisfies $||T_{yd}(z)||_\infty\le\gamma$ if and only if there exists a matrix $P\in\mathbb{S}^{n_x}_{++}$ such that
    \begin{equation}\label{Hinf only involved AB}
        \begin{bmatrix}\centering
	P		&(A+BK)P		&G		&0\\
	P(A+BK)^\top		&P		&0		&P(C+DK)^\top\\
	G^\top		&0		&I		&H^\top\\
    0          &(C+DK)P    &H  &\gamma^2 I
    \end{bmatrix}  \succ 0.
    \end{equation}
\end{mylem}

\subsection{Structured controllers}
The controller constructs a feedback loop that connects system states and control inputs. The optimal controller is typically implemented in a centralized manner, requiring each local controller to access full-state information. However, the centralized controller is not feasible for various applications because it can pose a significantly high computational burden on large-scale systems, which often have limited capabilities. This motivates the design of sparse controllers, where only partial state information is utilized to compute each local controller.

In this paper, we are interested in the design of controllers with \textit{a priori} structural constraints $K\in S$, where $S$ is a linear subspace. The corresponding structured identity is described by a binary matrix $I_S$ defined by 
\begin{equation}
    [I_S]_{ij} \triangleq \begin{cases}
            1, &\text{if} \ K_{ij}\ \text{is a free variable},\\
            0, &\text{if} \ K_{ij} = 0\ \text{is required}.
        \end{cases}
\end{equation}
Then we define the complementary structure identity matrix as $I_{S^c} \triangleq \mathbf{1} - I_S$. Consequently, the structural constraint can be expressed as $K\circ I_{S^c} = 0$. It is worth noting that, there are no general conditions guaranteeing the feasibility of a structured stabilizing controller, since such a controller may not exist.
\section{Data collection and problem formulation}\label{Data collection and problem formulation}
This section presents the data collection mechanism, which is fundamental for further analysis. Here we adopt the same setting as in \cite{miller2024data}. We assume the matrices $A, B$ are deterministic but unknown \textit{a priori}, and the matrices $C,D,G,H$ are deterministic and known \textit{a priori}. Besides, we can access a trajectory of data that we sample with known control inputs and unknown bounded noise.
\subsection{Data collection}
Consider a linear input-state model given by
\begin{equation}
    x_{k+1} = A_*x_k + B_*u_k + w_k,
\end{equation}
where $A_*$ and $B_*$ are unknown ground truth system matrices in \eqref{basic system}, and $w_k$ is unknown process noise. We assume we can only access measured data of states and inputs, which are collected into the following matrices
\begin{equation}
\begin{aligned}
    X  &\triangleq \begin{bmatrix}x_0 &x_1 &\dots &x_T\end{bmatrix} \in \mathbb{R}^{n\times (T+1)},\\
    U &\triangleq \begin{bmatrix}u_0 &u_1 &\dots &u_{T-1}\end{bmatrix} \in \mathbb{R}^{m\times T}.
\end{aligned}
\end{equation}
Additionally, we make the following assumption.
\begin{myass}\label{bounded noise}
    The process noise is bounded, i.e., $||w_k||_2\le \epsilon$ for $\forall k=0,1,\ldots,T-1$ for some $\epsilon>0.$
\end{myass}
There are typically many system matrices $(A,B)$ that can generate $X,U$ with $w_k$ satisfying \cref{bounded noise}. We refer to them as \textit{all possible $(A,B)$ consistent with data}. {\color{blue}Miller et al. \cite{miller2024data} originally considered the energy-bounded noise, and they claimed this framework can be used to overapproximate the instantaneous-bounded noise.} However, this approach is quite conservative as $T$ increases. Instead, we directly consider the instantaneous-bounded noise, which is first studied in \cite{bisoffi2021trade}. Based on the data $(x_i,u_i,x_{i+1})$ at the \textit{i}-th timestep, we can describe a set of $(A,B)$
\begin{equation}
    \Sigma_i\triangleq \left\{(A,B) |\ x_{i+1} = Ax_i+ Bu_i + w_i, w_iw_i^\top\preceq \epsilon^2 I \right\}.
\end{equation}
Namely, all $(A,B)\in\Sigma_i$ can possibly generate $(x_i,u_i,x_{i+1})$. We now focus on $w_iw_i^\top\preceq \epsilon^2 I$. After eliminating $w_i$, we get
\begin{equation}\label{AB data QMI}
    \begin{bmatrix}
        I\\
        A^\top\\
        B^\top
    \end{bmatrix}^\top
    \begin{bmatrix}
        I &x_{i+1}\\
        0 &-x_i\\
        0 &-u_i
    \end{bmatrix}
    \begin{bmatrix}
        \epsilon^2 I &0\\
        0 &-1
    \end{bmatrix}
    \begin{bmatrix}
        I &x_{i+1}\\
        0 &-x_i\\
        0 &-u_i
    \end{bmatrix}^\top
    \begin{bmatrix}
        I\\
        A^\top\\
        B^\top
    \end{bmatrix}\succeq 0.
\end{equation}
We notice that the middle term of \eqref{AB data QMI} is constant and to facilitate discussions, we denote
\begin{equation}
    \Psi_i \triangleq\begin{bmatrix}
        I &x_{i+1}\\
        0 &-x_i\\
        0 &-u_i
    \end{bmatrix}
    \begin{bmatrix}
        \epsilon^2 I &0\\
        0 &-1
    \end{bmatrix}
    \begin{bmatrix}
        I &x_{i+1}\\
        0 &-x_i\\
        0 &-u_i
    \end{bmatrix}^\top.
\end{equation}
Clearly, the set of $(A,B)$ consistent with all data points is thus given by
\begin{equation}\label{system consistent to data}
    \Sigma = \bigcap_{i=0}^{T-1}{\left\{ (A,B) \Bigg|     \begin{bmatrix}
        I\\
        A^\top\\
        B^\top
    \end{bmatrix}^\top
    \Psi_i
    \begin{bmatrix}
        I\\
        A^\top\\
        B^\top
    \end{bmatrix}\succeq 0\right\}.}
\end{equation}
Our set of $(A,B)$ consistent with the data is precisely the real one, which is a subset of that obtained in \cite{miller2024data}. This result was demonstrated in \cite{bisoffi2021trade}. A larger possible set of $(A,B)$ introduces additional scenarios for the controller to consider, which could potentially degrade system performance. Therefore, our expression for the set of $(A,B)$ consistent with the data is less conservative and can lead to enhanced performance.

Now we are ready to formulate our problems.
\subsection{Problem formulation}
In this paper, we are interested in the data-driven optimal structured state-feedback controller design under $H_2$ and $H_\infty$ performance. Since we cannot determine the ground truth system given a sequence of measured data, the highest pursuit is to design a structured controller that minimizes the maximal $H_2(H_\infty)$ norm for $\forall(A,B)\in \Sigma$. To make comparisons, we also provide feasible and efficient solutions in the model-based structured controller design case and the data-driven unstructured controller design case. The problems are formulated as follows
\begin{itemize}
    \item \textbf{Model-based Structured $H_\infty$ Controllers Design}: When $A$ and $B$ are known, design a state-feedback controller $K\in S$, where $S$ is a given structure pattern, while minimizing the $H_\infty$ norm.
    \item \textbf{Data-driven Unstructured $H_2(H_\infty)$ Controllers Design}: When $A$ and $B$ are unknown, using the collected data $X,U$, design a state-feedback controller $K$ to minimize $\gamma$ such that $||T_{yd}(z)||_2\le\gamma$ ($||T_{yd}(z)||_\infty\le\gamma$) holds for $\forall(A,B)\in\Sigma$.
    \item \textbf{Data-driven Structured $H_2(H_\infty)$ Controllers Design}: When $A$ and $B$ are unknown, using the collected data $X,U$, design a state-feedback controller $K\in S$, where $S$ is a given structure pattern, to minimize $\gamma$ such that $||T_{yd}(z)||_2\le\gamma$ ($||T_{yd}(z)||_\infty\le\gamma$) holds for $\forall(A,B)\in\Sigma$.
\end{itemize}
\begin{myrem}
    For the second and third problems, we directly design an unstructured (structured) controller using the collected data $X,U$ to enhance the system's $H_2(H_\infty)$ performance. Since we do not know the ground truth $A_*,B_*$, we describe a set $\Sigma$ of $A,B$ that can possibly generate $X,U$ under the noise condition $w_k$ satisfying \cref{bounded noise}. Subsequently, we minimize $\gamma$ such that $||T_{yd}(z)||_2\le\gamma$ ($||T_{yd}(z)||_\infty\le\gamma$) holds for all $(A,B)\in\Sigma$ by leveraging the matrix S-procedure. However, this technique only gives sufficient conditions and thus is conservative. Therefore, the $\gamma^*$ returned by our algorithms represents an upper bound of the maximal $H_2(H_\infty)$ norm. We refer to this upper bound as the \textit{$H_2(H_\infty)$ performance bound}.
\end{myrem}

% {\color{blue}\section{Stabilization control}\label{stabilization control}
% In this section, we aim to design the sparsest stabilizing controller $K$ when $A,B$ are unknown using collected data $X,U$. In other words, we want to find the sparsest controller $K$ such that there exists $P\in\mathbb{S}^{n_x}_{++}$ and 
% \begin{equation}
%     P-(A+BK)P(A+BK) \succ 0
% \end{equation}
% for $\forall(A,B)\in\Sigma$. 
% }
\section{$H_2$ control}\label{H2 control}
In this section, we deal with the $H_2$ control case. We will solve $3$ sub-problems in $3$ different subsections respectively.
\subsection{Model-based structured controller design}
The problem of structured $H_2$ controller design in the model-based setting has been thoroughly investigated in the literature. Notable approaches include LMI-based methods and gradient-based algorithms, as discussed in \cite{fardad2014design,lin2011augmented}. In this work, we adopt the LMI-based technique from \cite{fardad2014design} to address the model-based structured $H_2$ controller design problem.
\subsection{Data-driven unstructured controllers design}
In this subsection, we aim to design the optimal unstructured controller $K$ when $A,B$ are unknown using collected data $X,U$. One powerful tool that facilitates our analysis is the (lossy) matrix S-procedure \cite{bisoffi2021trade}, a tool commonly used in data-driven control. We notice that the only equation involving $A,B$ is \eqref{H2 only involved AB}. Therefore, we first need to do some reorganizations and rewrite \eqref{H2 only involved AB} into a QMI in $A,B$. By doing so, we can establish a connection between \eqref{H2 only involved AB} and the collected data by using the matrix S-procedure, which is a direct extension of the S-procedure to the matrix case.

Following the approach in \cite{miller2024data}, by introducing $\beta>0$, \eqref{H2 only involved AB} can be reformulated as
\begin{equation}\label{AB model QMI}
    \begin{bmatrix}
        I\\
        A^\top\\
        B^\top
    \end{bmatrix}^\top \begin{bmatrix}
        P-GG^\top-\beta I  &0\\
        0       &-\begin{bmatrix}
            I\\
            K
        \end{bmatrix}P\begin{bmatrix}
            I\\
            K
        \end{bmatrix}^\top
    \end{bmatrix}   \begin{bmatrix}
        I\\
        A^\top\\
        B^\top
    \end{bmatrix} \succeq 0.
\end{equation}
Recall the form of the set of all possible systems consistent with data \eqref{system consistent to data}. The matrix S-procedure provides a sufficient condition for $\Sigma$ to imply \eqref{AB model QMI}, meaning \eqref{AB model QMI} holds for $\forall (A,B)\in \Sigma$. Throughout this paper, we use a common Lyapunov matrix $P$ for data-driven control.
\begin{mypro}\label{h2 s procedure}
    Equation \eqref{H2 only involved AB} holds for $\forall (A,B)\in \Sigma$ if there exists nonnegative $\alpha_0,\alpha_1,\ldots,\alpha_{T-1}$ and $\beta>0$ such that the following equation holds
    \begin{equation}\label{h2 s procedure equation}
         \begin{bmatrix}
        P-GG^\top-\beta I  &0\\
        0       &-\begin{bmatrix}
            I\\
            K
        \end{bmatrix}P\begin{bmatrix}
            I\\
            K
        \end{bmatrix}^\top
    \end{bmatrix} - \sum_{i=0}^{T-1}\alpha_i\Psi_i \succeq0.
    \end{equation}
\end{mypro}
\quad \textit{Proof}: Pre- and post-multiplying \eqref{h2 s procedure equation} with $\begin{bmatrix}
        I &A &B
    \end{bmatrix}$ and $\begin{bmatrix}
        I &A &B
    \end{bmatrix}^\top$, we can show that \eqref{AB model QMI} and further \eqref{H2 only involved AB} hold using \eqref{h2 s procedure equation} and \eqref{AB data QMI}. The proof is complete.
$\hfill \square$

\cref{h2 s procedure} is essential for the data-driven $H_2$ control, as it ensures the feasibility of \eqref{H2 only involved AB} directly through data, avoiding discussing all possible $(A,B)$ consistent with data.

Recall that our objective is to minimize $\gamma$ such that $||T_{yd}(z)||_2\le\gamma$ holds for all $(A,B)\in\Sigma$. To achieve this, we minimize $\gamma$ subject to the constraints in \eqref{H2 feasiblity}, replacing \eqref{H2 only involved AB} with the matrix S-procedure condition \eqref{h2 s procedure equation}. By introducing the change of variables $L=KP$ and applying the Schur complement, we formulate the following SDP problem:
\begin{subequations}\label{H2 data-driven unstructured}
\begin{align}
&\mathop{\min}_{P,L,Q,\{\alpha_i\}_{i=0,1,\ldots,T-1},\beta,\gamma}   \gamma  \qquad \text{s.t.} \\
     &  \begin{bmatrix}
        \setlength{\dashlinegap}{0.8pt}
        \begin{array}{ccc:c}
        P-GG^\top-\beta I & 0 &0 &0\\
        0 &0 &0 &P\\
        0 &0 &0 &L\\
        \hdashline
        0 &P &L^\top &P
        \end{array}
        \end{bmatrix}-\begin{bmatrix}
        \setlength{\dashlinegap}{0.8pt}
        \begin{array}{c:c}
        \sum_{i=0}^{T-1}\alpha_i\Psi_i &0\\
        \hdashline
        0 &0
        \end{array}
\end{bmatrix} \succeq 0,\\
    & \begin{bmatrix}\centering
	Q		&CP+DL		\\
	(CP+DL)^\top		&P\\
  \end{bmatrix} \succ 0,\operatorname{Tr}(Q)\le \gamma^2,\\
    & \alpha_0 \ge 0, \alpha_1 \ge 0,\ldots,\alpha_{T-1} \ge 0,\beta > 0,\\
    & P\in\mathbb{S}^{n_x}_{++}, Q\in\mathbb{S}^{n_y}_{++}.
\end{align}
\end{subequations}

We design the data-driven unstructured $H_2$ controller by solving Problem \eqref{H2 data-driven unstructured}. We denote the solution to Problem~\eqref{H2 data-driven unstructured} to be $(P^*,L^*,Q^*,\{\alpha^*_i\}_{i=0,1,\ldots,T-1},\beta^*,\gamma^*)$. The optimal controller is computed as $K^*=L^*{P^*}^{-1}$. The $\gamma^*$ is the minimal $\gamma$ such that $||T_{yd}(z)||_2\le\gamma$ for $\forall (A,B)\in \Sigma$. This case is considered in the first case of \cref{Unstructured optimal controllers design algorithm}.

\begin{myrem}\label{H2 unstructured advantage}
To the best of our knowledge, van Waarde et al. \cite{van2020noisy} were among the first to address data-driven unstructured $H_2$ controller design {\color{blue}under energy-bounded noise, and they claimed this framework could be used to tackle instantaneous bounds as a special case of energy bounds}. This methodology was subsequently adopted in the unstructured case by \cite{miller2024data}. {\color{blue}However, such energy-based overapproximation} introduces additional conservatism compared to directly handling instantaneous-bounded noise via the matrix S-procedure, as demonstrated in \cite{bisoffi2021trade}. {\color{blue}Our approach that leverages the matrix S-procedure to directly address instantaneous-bounded noise is therefore advantageous compared to the methods in \cite{van2020noisy} and \cite{miller2024data}.} This advantage is established theoretically in \cref{comparison with miller} and validated numerically in \cref{Simulations}.
\end{myrem}

\subsection{Data-driven structured controllers design}
In this subsection, we focus on designing an optimal structured controller $K\in S$ when $A,B$ are unknown. When incorporating $K\in S$ into the constraints, the change of variable technique in Problem \eqref{H2 data-driven unstructured} is no longer viable since this will lead to an intractable $LP^{-1}\in S$ constraint. To address this issue, we reformulate Problem~\eqref{H2 data-driven unstructured} into the following problem, in which $K$ is expressed explicitly so that the structural constraint $K\in S$ can be incorporated.
\begin{subequations}\label{H2 data-driven nonconvex structured}
\begin{align}
&\mathop{\min}_{P,K,Q,Y,\{\alpha_i\}_{i=0,1,\ldots,T-1},\beta,\gamma}   \gamma  \qquad \text{s.t.} \\
     & 
     \begin{bmatrix}
     \scalebox{0.92}{$
        \setlength{\dashlinegap}{0.8pt}
        \begin{array}{ccc:c}
        P-GG^\top-\beta I & 0 &0 &0\\
        0 &0 &0 &I\\
        0 &0 &0 &K\\
        \hdashline
        0 &I &K^\top &Y
        \end{array}$
        }
        \end{bmatrix}
        -
        \scalebox{0.98}{$
        \begin{bmatrix}
        \setlength{\dashlinegap}{0.8pt}
        \begin{array}{c:c}
        \sum_{i=0}^{T-1}\alpha_i\Psi_i &0\\
        \hdashline
        0 &0
        \end{array}
\end{bmatrix}$
} 
\succeq 0, \label{H2 data-driven nonconvex structured con1}\\
    & \begin{bmatrix}\centering
	Q		&C+DK		\\
	(C+DK)^\top		&Y\\
  \end{bmatrix} \succ 0, \operatorname{Tr}(Q)\le \gamma^2, K\circ I_{S^c} = 0,\label{H2 data-driven nonconvex structured con2}\\
    & Y \preceq P^{-1}, \label{H2 data-driven nonconvex structured con3}\\
    & \alpha_0 \ge 0, \alpha_1 \ge 0,\ldots,\alpha_{T-1}\ge0, \beta > 0, \label{H2 data-driven nonconvex structured con4}\\
    & P\in\mathbb{S}^{n_x}_{++}, Q\in\mathbb{S}^{n_y}_{++}, Y\in\mathbb{S}^{n_x}_{++}. \label{H2 data-driven nonconvex structured con5}
\end{align}
\end{subequations}
The primary challenge in Problem \eqref{H2 data-driven nonconvex structured} lies in the constraint $Y \preceq P^{-1}$, which introduces non-convexity in the feasibility region w.r.t. $(Y,P)$. To address this, we employ linearizing $P^{-1}$ around a prescribed point $\Tilde{P}$ \cite{fardad2014design} to relax the constraint and reformulate the problem as an SDP. Specifically, for a given point $\Tilde{P}$, $Y \preceq P^{-1}$ can be relaxed by linearizing $P^{-1}$ around $\Tilde{P}$ :
\begin{equation}\label{relax by linearization}
    Y \preceq \Tilde{P}^{-1} - \Tilde{P}^{-1}(P-\Tilde{P})\Tilde{P}^{-1}.
\end{equation}
Thus, Problem~\eqref{H2 data-driven nonconvex structured} can be relaxed to be
\begin{subequations}\label{H2 data-driven reformulated convex}
\begin{align}
&\mathop{\min}_{P,K,Q,Y,\{\alpha_i\}_{i=0,1,\ldots,T-1},\beta,\gamma}   \gamma  \qquad \text{s.t.} \\
         & Y \preceq \Tilde{P}^{-1} - \Tilde{P}^{-1}(P-\Tilde{P})\Tilde{P}^{-1},\\
    & \eqref{H2 data-driven nonconvex structured con1}, \eqref{H2 data-driven nonconvex structured con2}, \eqref{H2 data-driven nonconvex structured con4}, \eqref{H2 data-driven nonconvex structured con5}.
\end{align}
\end{subequations}
Problem~\eqref{H2 data-driven reformulated convex} is an SDP problem and can be solved using commercial solvers. We now characterize the relationship between Problem~\eqref{H2 data-driven nonconvex structured} and Problem~\eqref{H2 data-driven reformulated convex}.

 \begin{mythr}\label{relationship}
    Any feasible solution $(P,K,Q,Y,\\ \{\alpha_i\}_{i=0,1,\ldots,T-1},\beta,\gamma)$ to Problem~\eqref{H2 data-driven reformulated convex} is also feasible for Problem~\eqref{H2 data-driven nonconvex structured}. Given a feasible solution $(P,K,Q,Y,\{\alpha_i\}_{i=0,1,\ldots,T-1},\beta,\gamma)$ to Problem~\eqref{H2 data-driven nonconvex structured}, take $\Tilde{P} = P$, then $(P,K,Q,Y,\{\alpha_i\}_{i=0,1,\ldots,T-1},\beta,\gamma)$ is also feasible for Problem~\eqref{H2 data-driven reformulated convex}.
\end{mythr}
\quad \textit{Proof}: {\color{blue}The proof of the first claim is equivalent to showing that $\Tilde{P}^{-1} - \Tilde{P}^{-1}(P-\Tilde{P})\Tilde{P}^{-1} \preceq P^{-1}$, which follows from the convexity of the inverse mapping in the Loewner
order \cite{bhatia2013matrix}. Define $f(P) = P^{-1}$ on the cone of positive definite matrices. The inverse mapping is twice continuously differentiable, and its second-order directional derivative at $P$ along a symmetric direction $H$ is
\begin{equation}
        D^2(P^{-1})[H,H]=2P^{-1}HP^{-1}HP^{-1}\succeq 0.
\end{equation}
Hence, for any $X,Y\succ0$, define
\begin{equation}
    g(t) = f(X+t(Y-X)), t\in [0,1].
\end{equation}
By the chain rule,
\begin{equation}
    g''(t) = D^2f(X+t(Y-X))[Y-X,Y-X]\succeq 0, \forall t\in [0,1].
\end{equation}
Therefore $g(t)$ is convex in the Loewner order. By the first-order characterization of convex functions,
\begin{equation}
    g(1) \succeq g(0) + g'(0),
\end{equation}
which yields
\begin{equation}
    f(Y) \succeq f(X) + Df(X)[Y-X].
\end{equation}
Taking $X = \Tilde{P}$ and $Y = P$ and noting that $Df(X)[H]=-X^{-1}HX^{-1}$, we obtain 
\begin{equation}
    \Tilde{P}^{-1} - \Tilde{P}^{-1}(P-\Tilde{P})\Tilde{P}^{-1} \preceq P^{-1}.
\end{equation}
The proof of the first claim is complete. For the second claim, if we take $\Tilde{P} = P$, the constraint $Y \preceq \Tilde{P}^{-1} - \Tilde{P}^{-1}(P-\Tilde{P})\Tilde{P}^{-1}$ reduces to $Y \preceq P^{-1}$, which is satisfied since $(P,K,Q,Y,\gamma)$ is feasible for Problem (17). Therefore, $(P,K,Q,Y,\gamma)$ is feasible for Problem (19), and the proof is concluded. }$\hfill \square$

\cref{relationship} indicates that the feasibility region of Problem~\eqref{H2 data-driven reformulated convex} is contained within that of Problem~\eqref{H2 data-driven nonconvex structured}, which might cause infeasibility when incorporating the constraint $K \circ I_{S^c} = 0$. To tackle this difficulty, we slightly modify Problem~\eqref{H2 data-driven reformulated convex} into the following SDP problem

\begin{subequations}\label{H2 data-driven structured}
\begin{align}
&\mathop{\min}_{P,K,Q,Y,Z,\{\alpha_i\}_{i=0,1,\ldots,T-1},\beta,\gamma}   \gamma +\lambda\operatorname{Tr}(Z) \qquad \text{s.t.} \\
& \eqref{H2 data-driven nonconvex structured con1},\eqref{H2 data-driven nonconvex structured con2},\eqref{H2 data-driven nonconvex structured con4},\\
    & Y \preceq \Tilde{P}^{-1} - \Tilde{P}^{-1}(P-\Tilde{P})\Tilde{P}^{-1} + Z,\\
    & P\in\mathbb{S}^{n_x}_{++}, Q\in\mathbb{S}^{n_y}_{++}, Y\in\mathbb{S}^{n_x}_{++}, Z\in\mathbb{S}^{n_x}_{+}.
\end{align}
\end{subequations}

Problem \eqref{H2 data-driven structured} introduces an auxiliary variable $Z\in\mathbb{S}^{n_x}_{+}$ to enhance feasibility. The data-driven structured $H_2$ controller is obtained via an iterative algorithm, where $\Tilde{P}$ is updated at each iteration using the previous solution. The penalty parameter $\lambda$ is gradually increased to drive $Z$ towards zero, thereby recovering feasibility for Problem \eqref{H2 data-driven nonconvex structured}. This iterative procedure dynamically refines the feasibility region and can yield controllers with improved performance. The complete algorithm is presented in the first case of \cref{Structured optimal controllers design algorithm}, with $\delta$ serving as a safeguard to prevent $\lambda$ from becoming excessively large.

\section{{\color{blue}Theoretical properties analysis}}\label{comparison with miller}
{\color{blue}In this section, we analyze several fundamental properties of the proposed data-driven control framework. In particular, we study: (i) the relationship between different data-driven control frameworks and their induced feasibility regions, (ii) the behavior of the performance bound w.r.t. data length, and (iii) the feasibility exploration of different structured control techniques.

We first consider (i). Miller et al. \cite{miller2024data} considered energy-bounded noise and verified their proposed framework using instantaneous-bounded noise. The problem they solved is a special case of Problem \eqref{H2 data-driven unstructured} with $\alpha_0=\alpha_1=\cdots=\alpha_{T-1}$, so $\{\alpha_i\}_{i=0,1,\ldots,T-1}$ can be replaced by a single variable $\alpha$. By denoting their problem of interest as Problem (M), we have the following result:}

\begin{mythr}\label{better performance}
    The optimal $\gamma^*$ returned by Problem \eqref{H2 data-driven unstructured} is less than or equal to that returned by Problem (M).
\end{mythr}
\quad \textit{Proof}: We only need to show that the feasibility region of Problem \eqref{H2 data-driven unstructured} contains that of Problem (M). Consider a feasible solution $(P,L,Q,\alpha,\beta,\gamma)$ to Problem (M). Clearly, the same $P,L,Q,\beta,\gamma$ and $\alpha_0,\alpha_1,\ldots,\alpha_{T-1}=\alpha$ provide a feasible solution to Problem \eqref{H2 data-driven unstructured}. The proof is complete.
$\hfill \square$

{\color{blue}Next, we consider (ii). We prove that in our proposed framework, incorporating additional data cannot deteriorate the performance bound, as described in the following theorem.}
\begin{mythr}\label{monotonicity}
    The performance bound $\gamma^*$ returned by Problem \eqref{H2 data-driven unstructured} is monotonically nonincreasing w.r.t. $T$.
\end{mythr}
\quad \textit{Proof}: Let the solution of Problem \eqref{H2 data-driven unstructured} be $(P^*,L^*,Q^*,\allowbreak\{\alpha^*_i\}_{i=0,1,\ldots,T-1},\beta^*,\gamma^*)$. Assume another data point is collected, increasing the data length to $T+1$. It follows that $(P^*,L^*,Q^*,\{\alpha^*_i\}_{i=0,1,\ldots,T-1},\alpha^*_T=0,\beta^*,\gamma^*)$ is feasible for the new problem. Thus, the optimal solution to the new problem is at least as good as the previous solution. The proof is complete. $\hfill \square$

{\color{blue}We then consider (iii). Classical structured control techniques typically leverage the results from \cite{ferrante2019design}. They utilize the change of variables technique and consider only a convex subset of the feasibility region for $K$ to satisfy a specified structure. Consequently, this method is conservative, leaving much of the feasibility region unexplored. In contrast, our methods can explore a larger feasibility region. By iteratively adjusting the considered convex feasibility region and improving the controller over time, our methods may lead to a structured controller with improved performance.

In summary, the proposed data-driven structured $H_2$ control framework integrates the advantages established in (i), (ii), and (iii). Specifically, the less conservative data-driven feasibility region in (i) enables improved performance bounds, the monotonicity property in (ii) guarantees that the performance bound does not deteriorate as more data is collected, and the enhanced feasibility exploration in (iii) allows for the design of structured controllers with superior performance. These theoretical benefits are further substantiated by our simulation results, which demonstrate that the proposed framework can achieve improved performance compared with existing literature.}

\section{$H_\infty$ control}\label{Hinf control}
In this section, we deal with the $H_\infty$ control case. Similar to the $H_2$ control case, we solve $3$ sub-problems in $3$ different subsections respectively.
\subsection{Model-based structured controllers design}
In this subsection, we assume $A,B$ are known. Similar to the $H_2$ case, the change of variable cannot be used. Therefore, we express $K$ explicitly to incorporate the structural constraint. Following the idea and technique in the $H_2$ data-driven structured controller design, we can formulate the following problem
\begin{subequations}\label{Hinf model-based structured}
\begin{align}
&\mathop{\min}_{P,K,Y,Z,\gamma}   \gamma +\lambda\operatorname{Tr}(Z) \qquad \text{s.t.} \\
    & \begin{bmatrix}\centering
	P		&A+BK		&G		&0\\
	(A+BK)^\top		&Y		&0	&(C+DK)^\top	\\
	G^\top		&0		&I		&H^\top\\
    0   &C+DK   &H  &\gamma^2 I
    \end{bmatrix}  \succ 0,\\
    & Y \preceq \Tilde{P}^{-1} - \Tilde{P}^{-1}(P-\Tilde{P})\Tilde{P}^{-1} + Z,K \circ I_{S^c} = 0,\\
    & P\in\mathbb{S}^{n_x}_{++}, Q\in\mathbb{S}^{n_y}_{++},Y\in\mathbb{S}^{n_x}_{++}.
\end{align}
\end{subequations}
The detailed algorithm is outlined in the second case of \cref{Structured optimal controllers design algorithm}.

For comparison purposes, we formulate a reference problem that minimizes $\gamma$ subject to the constraints in \eqref{Hinf only involved AB}. By introducing the change of variables $L = KP$ and restricting $P$ to be a diagonal matrix, we obtain the following SDP:
\begin{subequations}\label{Hinf model-based structured for comparisons}
\begin{align}
&\mathop{\min}_{P,L,\gamma} \gamma  \qquad \text{s.t.}\\
    & L\circ I_{S^c}=0,\quad P\in\mathbb{S}^{n_x}_{++},\quad P \circ (\mathbf{1}-I)=0,\\
        &\begin{bmatrix}
    P		&AP+BL     &G  &0		\\
    (AP+BL)^\top		&P     &0      &(CP+DL)^\top\\
        G^\top         &0     &I      &H^\top\\
        0      &CP+DL      &H      &\gamma^2 I
    \end{bmatrix} \succ 0.
\end{align}
\end{subequations}
This approach is commonly adopted in the literature, as post-multiplication by a diagonal matrix preserves the prescribed sparsity pattern. Consequently, the controller $K = LP^{-1}$ inherits the structure of $L$ and thus satisfies the structural constraint.

\subsection{Data-driven unstructured controllers design}
In this subsection, we are interested in the design of the optimal unstructured controller $K$ when $A,B$ are unknown. Our data-driven unstructured $H_\infty$ control method considers the $H\neq0$ case, which was not addressed in \cite{van2020noisy}. We will still utilize the matrix S-procedure to analyze, meaning we need to reformulate \eqref{Hinf only involved AB} into a QMI in $A,B$. After multiplying both sides of \eqref{Hinf only involved AB} by $\mathrm{blkdiag}(I,P^{-1},I,I)$ we get
\begin{equation}
    \begin{bmatrix}\centering
	P		&A+BK		&G		&0\\
	(A+BK)^\top		&P^{-1}		&0	&(C+DK)^\top	\\
	G^\top		&0		&I		&H^\top\\
    0   &C+DK   &H  &\gamma^2 I
    \end{bmatrix}  \succ 0.\\
\end{equation}
By applying the elementary block matrix operations, this equation can further be reformulated as 
\begin{equation}
    \begin{bmatrix}\centering
	P^{-1}		&0		&(A+BK)^\top		&(C+DK)^\top\\
	0		&I		&G^\top	&H^\top	\\
	A+BK		&G		&P		&0\\
    C+DK   &H   &0  &\gamma^2 I
    \end{bmatrix}  \succ 0,\\
\end{equation}
which is equivalent to the following equation through  the Schur complement
\begin{equation}
    \begin{bmatrix}
        P &0\\
        0 &\gamma^2 I
    \end{bmatrix}
    -\begin{bmatrix}
        A+BK &G\\
        C+DK &H
    \end{bmatrix}
    \begin{bmatrix}
        P &0\\
        0 &I
    \end{bmatrix}
    \begin{bmatrix}
        A+BK &G\\
        C+DK &H
    \end{bmatrix}^\top\succ 0.
\end{equation}
After some matrix product, we can put everything together
\begin{equation}\label{matrix products}
    \begin{bmatrix}
        P-A_K P A_K^\top-GG^\top &-A_K P C_K^\top-GH^\top\\
        -C_K P A_K^\top-HG^\top &\gamma^2 I - C_K P C_K^\top - HH^\top
    \end{bmatrix} \succ 0,
\end{equation}
where $A_K = A+BK, C_K = C+DK$. Two equations can be derived from \eqref{matrix products} 
\begin{subequations}
    \begin{align}
        &\begin{aligned}
         &P-A_K P A_K^\top-GG^\top - (A_K P C_K^\top+GH^\top)\\
         &\cdot(\gamma^2 I - C_K P C_K^\top - HH^\top)^{-1}(C_K P A_K^\top + HG^\top) \succ 0,
        \end{aligned} \label{useful}\\
        &\gamma^2 I - C_K P C_K^\top - HH^\top \succ 0.\label{useless}
    \end{align}
\end{subequations}
In the later analysis, we mainly focus on \eqref{useful} because \eqref{useless} is implied by \eqref{useful} in the later discussions. After introducing $\beta>0$, we can further formulate \eqref{useful} as
\begin{equation}
    \begin{bmatrix}
        I\\
        A^\top\\
        B^\top
    \end{bmatrix}^\top \Omega
    \begin{bmatrix}
        I\\
        A^\top\\
        B^\top
    \end{bmatrix} \succeq 0,
\end{equation}
% \begin{equation}
%     \begin{bmatrix}
%         I &A &B
%     \end{bmatrix} \Omega
%     \begin{bmatrix}
%         I &A &B
%     \end{bmatrix}^\top \succ 0
% \end{equation}
where 
\begin{equation*}
\scalebox{0.95}{$
\begin{aligned}
    &\Omega =  \begin{bmatrix}
        P-GG^\top-\beta I &0\\
        0 &-\begin{bmatrix}
            I\\
            K
        \end{bmatrix}P\begin{bmatrix}
            I\\
            K
        \end{bmatrix}^\top
    \end{bmatrix} - \\
    & \begin{bmatrix}
        GH^\top\\
        \begin{bmatrix}
            I\\
            K
        \end{bmatrix}PC_K^\top
    \end{bmatrix}
    (\gamma^2 I - C_K P C_K^\top - HH^\top)^{-1}
    \begin{bmatrix}
        GH^\top\\
        \begin{bmatrix}
            I\\
            K
        \end{bmatrix}PC_K^\top
    \end{bmatrix}^\top.
\end{aligned}
$}
\end{equation*}
The derivation is more involved than in the $H_2$ case, as it requires the congruence transformations and the elementary block matrix operations.
Similar to \cref{h2 s procedure}, \cref{hinf s procedure} follows from the matrix S-procedure.
\begin{mypro}\label{hinf s procedure}
    Equation \eqref{Hinf only involved AB} holds for $\forall (A,B)\in \Sigma$ if there exists nonnegative $\alpha_0,\alpha_1,\ldots,\alpha_{T-1}$ and $\beta>0$ such that the following equation holds
    \begin{equation}\label{hinf s procedure equation}
         \Omega - \sum_{i=0}^{T-1}\alpha_i\Psi_i \succeq0.
    \end{equation}
\end{mypro}
\quad \textit{Proof}: It is similar to the $H_2$ case and is omitted. $\hfill \square$

We first apply the Schur complement to \eqref{hinf s procedure equation}

\begin{equation}\label{hinf complex transformation}
\begin{aligned}
    &\begin{bmatrix}
    \setlength{\dashlinegap}{0.8pt}
    \begin{array}{c:c}
        \begin{bmatrix}
            \times &0\\
            0 &-\begin{bmatrix}
        I\\
        K
    \end{bmatrix}P\begin{bmatrix}
        I\\
        K
    \end{bmatrix}^\top
        \end{bmatrix} &\begin{bmatrix}
            \times\\
            \begin{bmatrix}
        I\\
        K
    \end{bmatrix}P\begin{bmatrix}
        I\\
        K
    \end{bmatrix}^\top \times
        \end{bmatrix}\\
        \hdashline\\[-13pt]
        \star & \times - \times\begin{bmatrix}
        I\\
        K
    \end{bmatrix}P\begin{bmatrix}
        I\\
        K
    \end{bmatrix}^\top \times
    \end{array}
    \end{bmatrix} \\
    &-\begin{bmatrix}
        \setlength{\dashlinegap}{0.8pt}
        \begin{array}{c:c}
        \sum_{i=0}^{T-1}\alpha_i\Psi_i &0\\
        \hdashline
        0 &0
        \end{array}
        \end{bmatrix} \succeq 0,
\end{aligned}
\end{equation}
where $\times$ represents constant and irrelevant matrix blocks, and $\star$ the symmetric block. Clearly, we can equivalently reformulate \eqref{hinf complex transformation} by introducing a new variable
\begin{equation}\label{introduce Y}
    R \succeq \begin{bmatrix}
            I\\
            K
        \end{bmatrix}P\begin{bmatrix}
            I\\
            K
        \end{bmatrix}^\top.
\end{equation}
% \begin{equation}\label{introduce Y}
%     R \succeq \begin{bmatrix}
%             I &K^\top
%         \end{bmatrix}^\top P\begin{bmatrix}
%             I &K^\top
%         \end{bmatrix}
% \end{equation}
After changing the variable $L=KP$ and minimizing $\gamma$, we obtain the following SDP problem
\begin{subequations}\label{Hinf data-driven unstructured}
\begin{align}
&\mathop{\min}_{P,L,R,\{\alpha_i\}_{i=0,1,\ldots,T-1},\beta,\gamma}   \gamma  \qquad \text{s.t.} \\
    &  \begin{aligned}
    &{\renewcommand{\arraystretch}{1.8}
        \scalebox{0.89}{$\begin{bmatrix}
        \setlength{\dashlinegap}{0.8pt}
        \begin{array}{cc:c}
        P-GG^\top-\beta I & 0 &GH^\top \\
        0 &-R &R{\renewcommand{\arraystretch}{1.5}\begin{bmatrix}
            C^\top\\
            D^\top
        \end{bmatrix}} \\
        \hdashline
        &\\[-16pt]
        HG^\top &{\renewcommand{\arraystretch}{1.5}\begin{bmatrix}
            C &D
        \end{bmatrix}}R &\gamma^2 I - 
        {\renewcommand{\arraystretch}{1.5}\begin{bmatrix}C^\top  \\D^\top\end{bmatrix}^\top} R
        {\renewcommand{\arraystretch}{1.5}\begin{bmatrix}C^\top  \\D^\top\end{bmatrix}} 
        - HH^\top
        \end{array}
        \end{bmatrix}$} }\\
        &-\begin{bmatrix}
        \setlength{\dashlinegap}{0.8pt}
        \begin{array}{c:c}
        \sum_{i=0}^{T-1}\alpha_i\Psi_i &0\\
        \hdashline
        0 &0
        \end{array}
        \end{bmatrix} \succeq 0,
        \end{aligned} \label{Hinf data-driven unstructured con1}\\
        &\begin{bmatrix}
        R - \begin{bmatrix}
            P &L^\top\\
            L &0
        \end{bmatrix} &\begin{bmatrix}
            0\\
            L
        \end{bmatrix}\\
        \begin{bmatrix}
            0 &L^\top
        \end{bmatrix} &P
    \end{bmatrix} \succ 0,\label{Hinf data-driven unstructured con2}\\
    & \alpha_0 \ge 0, \alpha_1 \ge 0,\ldots,\alpha_{T-1}\ge0, \beta > 0,\label{Hinf data-driven unstructured con3}\\
    & P\in\mathbb{S}^{n_x}_{++}, R\in\mathbb{S}^{n_x+n_u}_{++}.\label{Hinf data-driven unstructured con4}
\end{align}
\end{subequations}

We design the data-driven unstructured $H_\infty$ controller by solving Problem \eqref{Hinf data-driven unstructured}. We denote the solution to Problem~\eqref{Hinf data-driven unstructured} to be $(P^*,L^*,R^*,\{\alpha^*_i\}_{i=0,1,\ldots,T-1},\beta^*,\gamma^*)$. The optimal controller is computed as $K^*=L^*{P^*}^{-1}$. The $\gamma^*$ is the minimal $\gamma$ such that $||T_{yd}(z)||_\infty\le\gamma$ for $\forall (A,B)\in \Sigma$. This case is incorporated as the second case of \cref{Unstructured optimal controllers design algorithm}.

For comparison, we extend the set characterization of all possible system matrices consistent with the data, as well as the methodology for handling all possibilities in \cite{miller2024data}, to the $H_\infty$ case. Specifically, the reference method corresponds to Problem \eqref{Hinf data-driven unstructured} with $\alpha_0 = \alpha_1 = \cdots = \alpha_{T-1}$, where the multipliers $\{\alpha_i\}_{i=0,1,\ldots,T-1}$ are replaced by a single variable $\alpha$. We denote this formulation as Problem (N).

Similar to the $H_2$ case, the performance bound returned by Problem \eqref{Hinf data-driven unstructured} outperforms that returned by Problem (N) and also behaves monotonicity property in the $H_\infty$ case, which is summarized in the following theorems.
\begin{mythr}
    The optimal $\gamma^*$ returned by Problem \eqref{Hinf data-driven unstructured} is less than or equal to that returned by Problem (N). 
\end{mythr}
\quad \textit{Proof}: It is similar to the $H_2$ case and is omitted. $\hfill \square$
\begin{mythr}\label{hinf monotonicity}
    The performance bound $\gamma^*$ returned by Problem \eqref{Hinf data-driven unstructured} is monotonically nonincreasing w.r.t. $T$.
\end{mythr}
\quad \textit{Proof}: It is similar to the $H_2$ case and is omitted. $\hfill \square$
\begin{myrem}\label{Hinf unstructured advantage}
    Compared to \cite{van2020noisy}, our data-driven unstructured $H_\infty$ control approach leverages the matrix S-procedure to handle instantaneous-bounded noise, resulting in reduced conservatism as discussed in \cref{H2 unstructured advantage}.
\end{myrem}

\subsection{Data-driven structured controllers design}
To extend the data-driven control from the unstructured case to the structured case, we only need to apply the Schur complement to \eqref{introduce Y} and linearize the inverse term. Following the idea in the $H_2$ control case, we directly provide the problem of interest below for conciseness. 
\begin{subequations}\label{Hinf data-driven structured}
\begin{align}
&\mathop{\min}_{P,R,Y,K,Z,\{\alpha_i\}_{i=0,1,\ldots,T-1},\beta,\gamma}   \gamma + \lambda\operatorname{Tr}(Z) \qquad \text{s.t.} \\
        & \eqref{Hinf data-driven unstructured con1}, \eqref{Hinf data-driven unstructured con3}\\
    & \begin{bmatrix}
        R  &\begin{bmatrix}
            I\\
            K
        \end{bmatrix}\\
        \begin{bmatrix}
            I &K^\top
        \end{bmatrix} &Y
    \end{bmatrix} \succ 0,K\circ I_{S^c} = 0\\
    & Y \preceq \Tilde{P}^{-1} - \Tilde{P}^{-1}(P-\Tilde{P})\Tilde{P}^{-1} + Z,\\
    & P\in\mathbb{S}^{n_x}_{++}, R\in\mathbb{S}^{n_x+n_u}_{++},Y\in\mathbb{S}^{n_x}_{++},Z\in\mathbb{S}^{n_x}_{+}.
\end{align}
\end{subequations}
The specifics of the algorithm are presented in the third case of \cref{Structured optimal controllers design algorithm}.

Similar to the previous cases, we construct a reference problem by combining Problem \eqref{Hinf model-based structured for comparisons} and Problem (N) to facilitate comparison with Problem \eqref{Hinf data-driven structured}. We denote this combined formulation as Problem (N'). In Problem (N'), we set $\alpha_0 = \alpha_1 = \cdots = \alpha_{T-1}$ and introduce a single multiplier $\alpha$ in place of $\{\alpha_i\}_{i=0,1,\ldots,T-1}$. Furthermore, we impose the constraint that $P$ is diagonal and enforce $L \circ I_{S^c} = 0$ to ensure the prescribed structural pattern.
\begin{algorithm}[t]
\caption{Unstructured optimal controllers design.}
  \begin{algorithmic}[1]
    \State \textbf{Output:} $K^*,\gamma^*$;
    \State $\text{(P)}=\begin{cases}
    \text{Problem \eqref{H2 data-driven unstructured}:}\\ 
            \qquad \quad \text{Data-driven}\ H_2\ \text{controllers design}\\
     \text{Problem \eqref{Hinf data-driven unstructured}:} \\
            \qquad \quad \text{Data-driven}\ H_\infty\ \text{controllers design}\\
    \end{cases}$
    
    \State Solve (P) and obtain the solution $L^*, P^*, \gamma^*$;
    \State \Return $L^* P^{*-1}, \gamma^*$. \label{return}
  \end{algorithmic}
  \label{Unstructured optimal controllers design algorithm}
\end{algorithm}

\begin{algorithm}[t]
\caption{Structured optimal controllers design.}
  \begin{algorithmic}[1]
    \State \textbf{Output:} $K^*,\gamma^*$;
    \State $\text{(P)}=\begin{cases}
    \text{Problem \eqref{H2 data-driven structured}:} \\
            \qquad \quad \text{Data-driven}\ H_2\ \text{controllers design}\\
    \text{Problem \eqref{Hinf model-based structured}:} \\
            \qquad \quad \text{Model-based}\ H_\infty\ \text{controllers design}\\
     \text{Problem \eqref{Hinf data-driven structured}:} \\
            \qquad \quad \text{Data-driven}\ H_\infty\ \text{controllers design}\\
    \end{cases}$
    \State Initialize $\lambda=\lambda_0, \mu>1, \delta>0, \epsilon>0, P_0 = I, k=0$;
    \Repeat
    % \vspace{2pt}
    \State Solve (P) with $\Tilde{P} = P_k$ and $\lambda$;
    \State Assign the solutions to $K_{k+1}, P_{k+1}, \gamma_{k+1}$;
    \State \textbf{if} $\lambda<\delta$ \textbf{then} $\lambda = \mu\lambda$ \textbf{end if}
    \State $k = k + 1$;
    % \vspace{-5pt}
    \Until{$||P_{k}-P_{k-1}||_F<\epsilon$ and $||Y_k-P_k^{-1}||_F<\epsilon$;\vspace{5pt}}
    \State \Return $K_k,\gamma_k$.
  \end{algorithmic}
\label{Structured optimal controllers design algorithm}
\end{algorithm}

\begin{myrem}
    In this paper, we do not claim the convergence of the iterative algorithm, as there may not be a stabilizing controller with the specified sparsity pattern, leading to potential infeasibility. However, given reasonable structure constraints, our extensive numerical examples consistently demonstrate that our algorithms yield a stabilizing structured controller with satisfactory performance. We believe that our methods hold promise for data-driven structured controller design.
\end{myrem}

\begin{myrem}\label{disadvantage}
    The primary limitation of our proposed methods lies in their computational complexity, since the algorithms introduce additional optimization variables and require iterative solution procedures. Furthermore, the number of variables grows with the length of the data sequence, leading to increased computational demands for larger datasets. Nevertheless, this limitation is offset by the substantial performance improvements demonstrated over existing approaches such as \cite{miller2024data}, as evidenced by our simulation results.
\end{myrem}

% \begin{table}[!t]
% \centering
% \caption{$H_2$ performance bound comparisons for $T=20$}
% \begin{tabular}{lccccr}
% \toprule
% Design & $(A,B)$ & $\epsilon=0.05$ & $\epsilon=0.1$  &$\epsilon=0.2$\\
% \midrule
% Unstructured in \cite{miller2024data}& 2.1537 & 2.3448 & 3.0939  &5.2814\\
% \textbf{Our Unstructured} &2.1537 &2.2663 &2.5660 &2.7889\\
% Structured in \cite{miller2024data} & 2,9794 & 3.5494 & 4.6806  &11.3186\\
% \textbf{Our Structured} & 2.7165 & 2.9154 & 3.2249  & 4.0422\\
% \bottomrule
% \end{tabular}
% \label{H2 comparison for fixed T}
% \end{table}

\begin{table}[!t]
\centering
\caption{$H_2$ performance bound comparisons for $T=20$.}
\scalebox{0.89}{
\begin{tabular}{lccccr}
\toprule
\multicolumn{2}{c}{Design}              & $(A,B)$ & $\epsilon=0.05$ & $\epsilon=0.1$  &$\epsilon=0.2$\\
\midrule
\multirow{2}{*}{Unstructured} & \cite{miller2024data} & 2.1537 & 2.3448 & 3.0939  &5.2814\\
                             & \textbf{Ours}      &2.1537 &2.2663 &2.5660 &2.7889\\
\hline
\multirow{4}{*}{Structured}  & \cite{miller2024data}  & 2.9794 & 3.5494 & 4.6806  &11.3186\\
&{\color{blue}\textbf{Ours}$^{\#}$}     & {\color{blue}2.7165} & {\color{blue}3.0760} & {\color{blue}3.7717}  & {\color{blue}6.1391}\\
& \cite{miller2024data}$^*$  & 2.9794 & 3.2442 & 3.6890  &5.4504\\
& \textbf{Ours}     & 2.7165 & 2.9154 & 3.2249  & 4.0422\\
\bottomrule
\end{tabular}}
\label{H2 comparison for fixed T}
\end{table}

\section{Simulations}\label{Simulations}
In this section, we will verify the effectiveness of our methods via several numerical examples related to $H_2$ control and $H_\infty$ control. {\color{blue}The code is available at https://github.com/Zhaohua0418/Data-Driven-Structured}. Specifically, we will compare our methods with those proposed in \cite{miller2024data} and show our improvement. The parameters mentioned in \cref{Structured optimal controllers design algorithm} are taken as $\lambda_0=1,\mu=2,\delta=1e8,\epsilon=0.01$. In the result tables, column $(A,B)$ represents the model-based case, and other columns represent the data-driven case. The values recorded are the $H_2(H_\infty)$ performance for the model-based case and the $H_2(H_\infty)$ performance bound for the data-driven case.

\begin{table}[!t]
\centering
\caption{$H_2$ performance bound comparisons for $\epsilon=0.1$.}
\scalebox{0.92}{
\begin{tabular}{lccccr}
\toprule
\multicolumn{2}{c}{Design}              & $(A,B)$ & $T = 6$ & $T= 10$ &$T = 15$\\
\midrule
\multirow{2}{*}{Unstructured} & \cite{miller2024data} & 2.1537 & 2.9911 & 2.8156  &2.6836\\
                             & \textbf{Ours}      &2.1537 &2.7494&2.5645& 2.4374\\
\hline
\multirow{4}{*}{Structured}  & \cite{miller2024data}  & 2.9794 & 4.4036 & 4.4323 &4.4456\\
& {\color{blue}\textbf{Ours}$^{\#}$}     &{\color{blue}2.7165} & {\color{blue}3.7204} & {\color{blue}3.5982}  &{\color{blue}3.5101}\\
                            & \cite{miller2024data}$^*$  & 2.9794 & 4.0991 & 3.6666 &3.5832\\
                                                         & \textbf{Ours}     &2.7165 & 3.5249 & 3.2429  &3.1485\\
\bottomrule
\end{tabular}}
\label{H2 comparison for fixed epsilon}
\end{table}

% \begin{table}[!t]
% \centering
% \caption{$H_2$ performance bound comparisons for $\epsilon=0.1$}
% \begin{tabular}{lccccr}
% \toprule
% Design & $(A,B)$ & $T = 6$ & $T= 10$ &$T = 15$\\
% \midrule
% Unstructured in \cite{miller2024data} & 2.1537 & 2.9911 & 2.8156  &2.6836\\
% \textbf{Our Unstructured} &2.1537 &2.7494&2.5645& 2.4374\\
% Structured in \cite{miller2024data} & 2,9794 & 4.4036 & 4.4323 &4.4456\\
% \textbf{Our Structured} & 2.7165 & 3.5249 & 3.2429  &3.1485\\
% \bottomrule
% \end{tabular}
% \label{H2 comparison for fixed epsilon}
% \end{table}

% \begin{table}[!t]
% \centering
% \caption{$H_\infty$ performance bound comparisons for $T=50$}
% \begin{tabular}{lccccr}
% \toprule
% Design & $(A,B)$ &$\epsilon=0.01$ & $\epsilon=0.05$  &$\epsilon=0.15$ \\
% \midrule
% Unstructured in \cite{miller2024data}& 0.7815 & 0.8035 & 0.9063  &1.6483\\
% \textbf{Our unstructured} &0.7815 &0.7921 &0.8207 &1.0303 \\
% $P$ is diagonal & 3.7464 & 4.0950 & 6.6234  & Infeasible \\
% \textbf{Our method} & 1.0580 & 1.0890 & 1.1826  &1.5969\\
% \bottomrule
% \end{tabular}
% \label{Hinf comparison for fixed T}
% \end{table}

\begin{table}[!t]
\centering
\caption{$H_\infty$ performance bound comparisons for $T=50$.}
\scalebox{0.84}{
\begin{tabular}{lccccr}
\toprule
\multicolumn{2}{c}{Design}              & $(A,B)$ &$\epsilon=0.01$ & $\epsilon=0.05$  &$\epsilon=0.15$\\
\midrule
\multirow{2}{*}{Unstructured} & \cite{miller2024data} & 0.7815 & 0.8035 & 0.9063  &1.6483\\
                             & \textbf{Ours}      &0.7815 &0.7921 &0.8207 &1.0303\\
\hline
\multirow{4}{*}{Structured}  & $P$ diag  & 3.7464 & 4.0950 & 6.6234  & Infeasible\\
& {\color{blue}\textbf{Ours}$^{\#}$}     & {\color{blue}1.0580} & {\color{blue}1.1323} & {\color{blue}1.4401}  &{\color{blue}4.1559}\\
                            & $P$ diag$^*$  & 3.7464 & 3.9149 & 4.7155  & 7.9693\\
                                                         & \textbf{Ours}     & 1.0580 & 1.0896 & 1.1825  &1.5969\\
                            \bottomrule
\end{tabular}}
\label{Hinf comparison for fixed T}
\end{table}

% \begin{table}[!t]
% \centering
% \caption{$H_\infty$ performance bound comparisons for $\epsilon=0.05$}
% \begin{tabular}{lccccr}
% \toprule
% Design & $(A,B)$ & $T = 10$ & $T= 20$ &$T = 40$\\
% \midrule
% Unstructured in \cite{miller2024data}& 0.7815 & 0.8802 & 0.9298  &0.9082\\
% \textbf{Our unstructured} &0.7815 &0.8521 &0.8306&0.8207 \\
% $P$ is diagonal & 3.7464 & 7.6607 & 7.8820 & 6.3521\\
% \textbf{Our method} & 1.0580 & 1.3903 & 1.2723  &1.1866\\
% \bottomrule
% \end{tabular}
% \label{Hinf comparison for fixed epsilon}
% \end{table}

\begin{table}[!t]
\centering
\caption{$H_\infty$ performance bound comparisons for $\epsilon=0.05$.}
\scalebox{0.9}{
\begin{tabular}{lccccr}
\toprule
\multicolumn{2}{c}{Design}              & $(A,B)$ & $T = 10$ & $T= 20$ &$T = 40$\\
\midrule
\multirow{2}{*}{Unstructured} & \cite{miller2024data} & 0.7815 & 0.8802 & 0.9298  &0.9082\\
                             & \textbf{Ours}      &0.7815 &0.8521 &0.8367&0.8207 \\
\hline
\multirow{3}{*}{Structured}  & $P$ diag & 3.7464 & 7.6607 & 7.8820 & 6.3521\\
& {\color{blue}\textbf{Ours}$^{\#}$}     & {\color{blue}1.0580} & {\color{blue}1.5233} & {\color{blue}1.5667}  &{\color{blue}1.4213}\\
                            & $P$ diag$^*$ & 3.7464 & 6.2570 & 4.9125 & 4.7494\\  
                                                         & \textbf{Ours}     & 1.0580 & 1.3911 & 1.2720  &1.1866\\
                            \bottomrule
\end{tabular}}
\label{Hinf comparison for fixed epsilon}
\end{table}
\begin{figure}[t]
\centering
\includegraphics[width=0.5\textwidth]{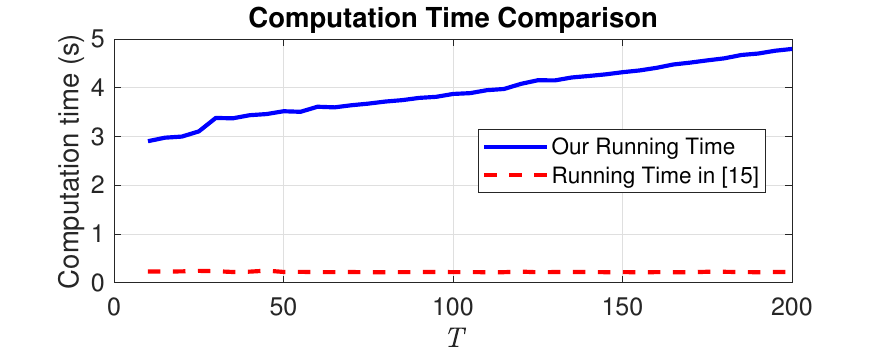}
\caption{Computation time comparison between our structured $H_2$ control and \cite{miller2024data} for different $T$ when $\epsilon=0.1$.}
\label{Computation time comparison}
\end{figure}

\subsection{$H_2$ control case study}
Consider a discrete-time LTI system \eqref{basic system} with system matrices provided in \cite{miller2024data}
\begin{equation*}
\begin{split}    
    &A = \begin{bmatrix}
        -0.4095 &0.4036 &-0.0874\\
        0.5154 &-0.0815 &0.1069\\
        1.6715 &0.7718 &-0.3376
    \end{bmatrix}, C = \begin{bmatrix}
        I_3\\
        0_{2\times3}
    \end{bmatrix},G=I_3,\\
    &B = \begin{bmatrix}
        0 &-0.6359 &-0.0325\\
        0 &-0.1098 &2.2795
    \end{bmatrix}^\top,D = \begin{bmatrix}
        0_{2\times3} &I_2
    \end{bmatrix}^\top,H=0.
\end{split}
\end{equation*}
We specify a structural constraint $K\in S$, with

\begin{equation}
    I_S = \begin{bmatrix}
        1 &1 &0\\
        0 &1 &1
    \end{bmatrix}.
\end{equation}

\cref{H2 comparison for fixed T} presents the $H_2$ performance comparisons for different design settings with $T=20$. {\color{blue}The notation``\cite{miller2024data}'' represents the approach proposed in \cite{miller2024data}. ``\textbf{Ours}$^\#$'' denotes testing our linearization approach within the energy-based overapproximation framework, as adopted in \cite{miller2024data}. In addition,} ``\cite{miller2024data}$^*$'' refers to the use of the matrix S-procedure to directly address instantaneous-bounded noise in \cite{miller2024data}, rather than employing overapproximation. {\color{blue}Finally, ``\textbf{Ours}'' refers to the method proposed in this paper.} As $\epsilon$ increases, the $H_2$ performance upper bound also increases, since greater uncertainty requires consideration of a larger set $\Sigma$. The comparisons between \cite{miller2024data} and ``\textbf{Ours}'' in the unstructured case, as well as between \cite{miller2024data} and \cite{miller2024data}$^*$ in the structured case, demonstrate the benefit of directly handling instantaneous-bounded noise, thereby validating our theoretical results in \cref{better performance}. {\color{blue}Furthermore, since the data-driven frameworks in \cite{miller2024data} and ``\textbf{Ours}$^{\#}$'', as well as those in \cite{miller2024data}$^*$ and ``\textbf{Ours}'', are identical,} the results indicate that the use of the linearization technique enables our ILMI algorithms to achieve controllers with improved performance. Additionally, \cref{H2 comparison for fixed epsilon} provides results for $\epsilon=0.1$, further confirming the effectiveness of our approach. It is noteworthy that the performance bound in \cite{miller2024data} does not decrease monotonically w.r.t. $T$, as shown in \cref{H2 comparison for fixed epsilon}. In contrast, our methods guarantee a monotonic decrease w.r.t. $T$, consistent with the findings in \cref{monotonicity}. However, the improved performance is accompanied by increased computational cost, as illustrated in \cref{Computation time comparison}. For $\epsilon=0.1$, we compare the computation time of our structured $H_2$ control method with that of \cite{miller2024data} for various values of $T$ ranging from $10$ to $200$. Our method requires more computation time due to the iterative procedure and the larger number of optimization variables. Moreover, the computation time increases with $T$, which is consistent with the discussion in \cref{disadvantage}.

\subsection{$H_\infty$ control case study}
Consider a discrete-time LTI system \eqref{basic system} with the following system matrices
\begin{equation*}
    \begin{split}
        &A = \begin{bmatrix}
            0.8 &0.2  &0.1\\
          0.1 &0.7 &-0.3\\
         -0.3 &0.5  &0.9
        \end{bmatrix},\ B = \begin{bmatrix}
            1 &0\\
            0 &1\\
            1 &1
        \end{bmatrix},\ G = \begin{bmatrix}
            0.3 &0.1\\
            0.2 &0.2\\
            0.1 &0.3
        \end{bmatrix},\\
        &C = I_3,\ D=\begin{bmatrix}
            0.1 &0.3 &0.2\\
            0.2 &0.1 &0.1
        \end{bmatrix}^\top,\ H = \begin{bmatrix}
            0.1 &0.2 &0.3\\
            0.1 &0.2 &0.3\
        \end{bmatrix}^\top.
    \end{split}
\end{equation*}
We specify a structural constraint $K\in S$, with
\begin{equation}
    I_S = \begin{bmatrix}
        1 &1 &0\\
        1 &1 &0
    \end{bmatrix}.
\end{equation}
The label ``\cite{miller2024data}'' refers to the extension of the data-driven unstructured controller design from \cite{miller2024data} to the $H_\infty$ case. ``$P$ diag'' denotes the approach where the data-driven framework from \cite{miller2024data} is adopted, with $P$ constrained to be diagonal and $L$ following the prescribed sparsity pattern. {\color{blue}Additionally, ``\textbf{Ours}$^{\#}$'' indicates our relaxation approach tested within the data-driven framework of \cite{miller2024data}.} The label ``$P$ diag$^*$'' refers to our proposed data-driven framework, using the same structural constraints as in ``$P$ diag''. {\color{blue}``\textbf{Ours}'' refers to our method.} The results in \cref{Hinf comparison for fixed T} and \cref{Hinf comparison for fixed epsilon} demonstrate that our methods consistently achieve improved $H_\infty$ performance bounds compared to the baselines. The underlying reasons are similar to those discussed for the $H_2$ case.

\section{Conclusion}\label{Conclusion}
In this paper, we focused on the data-driven optimal structured controller design, considering both $H_2$ performance and $H_\infty$ performance. For each criterion, we consider the model-based and data-driven, structured control and unstructured control, and developed promising algorithms for each case. {\color{blue} We demonstrate that our methods can achieved improved performance over those in the existing literature at the cost of computation via several numerical simulations.}

\footnotesize
\bibliographystyle{unsrtnat} % 数字引用格式  
\bibliography{reference}

@article{miller2024data,
  title={Data-Driven Structured Robust Control of Linear Systems},
  author={Miller, Jared and Eising, Jaap and D{\"o}rfler, Florian and Smith, Roy S},
  journal={arXiv preprint arXiv:2411.11542},
  year={2024}
}

@article{van2020noisy,
  title={From noisy data to feedback controllers: Nonconservative design via a matrix {S}-lemma},
  author={van Waarde, Henk J and Camlibel, M Kanat and Mesbahi, Mehran},
  journal={IEEE Transactions on Automatic Control},
  volume={67},
  number={1},
  pages={162--175},
  year={2020},
  publisher={IEEE}
}

@inproceedings{eising2022informativity,
  title={Informativity for centralized design of distributed controllers for networked systems},
  author={Eising, Jaap and Cort{\'e}s, Jorge},
  booktitle={Proceedings of European Control Conference (ECC)},
  pages={681--686},
  year={2022},
  organization={IEEE}
}

@article{bisoffi2022data,
  title={Data-driven control via {P}etersen’s lemma},
  author={Bisoffi, Andrea and De Persis, Claudio and Tesi, Pietro},
  journal={Automatica},
  volume={145},
  pages={110537},
  year={2022},
  publisher={Elsevier}
}

@article{bisoffi2021trade,
  title={Trade-offs in learning controllers from noisy data},
  author={Bisoffi, Andrea and De Persis, Claudio and Tesi, Pietro},
  journal={Systems \& Control Letters},
  volume={154},
  pages={104985},
  year={2021},
  publisher={Elsevier}
}

@article{blondel1997np,
  title={{NP}-hardness of some linear control design problems},
  author={Blondel, Vincent and Tsitsiklis, John N},
  journal={SIAM {J}ournal on {C}ontrol and {O}ptimization},
  volume={35},
  number={6},
  pages={2118--2127},
  year={1997},
  publisher={SIAM}
}

@inproceedings{polyak2013lmi,
  title={An {LMI} approach to structured sparse feedback design in linear control systems},
  author={Polyak, Boris and Khlebnikov, Mikhail and Shcherbakov, Pavel},
  booktitle={Proceedings of European Control Conference (ECC)},
  pages={833--838},
  year={2013},
  organization={IEEE}
}

@article{lin2011augmented,
  title={Augmented {L}agrangian approach to design of structured optimal state feedback gains},
  author={Lin, Fu and Fardad, Makan and Jovanovic, Mihailo R},
  journal={IEEE Transactions on Automatic Control},
  volume={56},
  number={12},
  pages={2923--2929},
  year={2011},
  publisher={IEEE}
}

@article{jovanovic2016controller,
  title={Controller architectures: Tradeoffs between performance and structure},
  author={Jovanovi{\'c}, Mihailo R and Dhingra, Neil K},
  journal={European {J}ournal of {C}ontrol},
  volume={30},
  pages={76--91},
  year={2016},
  publisher={Elsevier}
}

@article{yang2024log,
  title={Log-{B}arrier Search for Structural Linear Quadratic Regulators},
  author={Yang, Nachuan and Tang, Jiawei and Li, Yuzhe and Shi, Guodong and Shi, Ling},
  journal={IEEE Transactions on Automatic Control},
  volume={70},
  number={3},
  pages={1965-1972},
  year={2024},
  publisher={IEEE}
}

@inproceedings{fardad2014design,
  title={On the design of optimal structured and sparse feedback gains via sequential convex programming},
  author={Fardad, Makan and Jovanovi{\'c}, Mihailo R},
  booktitle={Proceedings of American Control Conference},
  pages={2426--2431},
  year={2014},
  organization={IEEE}
}

@article{ferrante2019design,
  title={On the design of structured stabilizers for {LTI} systems},
  author={Ferrante, Francesco and Dabbene, Fabrizio and Ravazzi, Chiara},
  journal={IEEE Control Systems Letters},
  volume={4},
  number={2},
  pages={289--294},
  year={2019},
  publisher={IEEE}
}

@article{de2019formulas,
  title={Formulas for data-driven control: Stabilization, optimality, and robustness},
  author={De Persis, Claudio and Tesi, Pietro},
  journal={IEEE Transactions on Automatic Control},
  volume={65},
  number={3},
  pages={909--924},
  year={2019},
  publisher={IEEE}
}

@inproceedings{berberich2020robust,
  title={Robust data-driven state-feedback design},
  author={Berberich, Julian and Koch, Anne and Scherer, Carsten W and Allg{\"o}wer, Frank},
  booktitle={Proceedings of American Control Conference (ACC)},
  pages={1532--1538},
  year={2020},
  organization={IEEE}
}

@article{van2020data,
  title={Data informativity: A new perspective on data-driven analysis and control},
  author={van Waarde, Henk J and Eising, Jaap and Trentelman, Harry L and Camlibel, M Kanat},
  journal={IEEE Transactions on Automatic Control},
  volume={65},
  number={11},
  pages={4753--4768},
  year={2020},
  publisher={IEEE}
}

@book{zhou1996robust,
  title={Robust and Optimal Control},
  author={Zhou, K. and Doyle, J.C. and Glover, K.},
  year={1996},
  publisher={Prentice Hall}
}

@article{gahinet1994linear,
  title={A linear matrix inequality approach to {$H_\infty$} control},
  author={Gahinet, Pascal and Apkarian, Pierre},
  journal={International {J}ournal of {R}obust and {N}onlinear {C}ontrol},
  volume={4},
  number={4},
  pages={421--448},
  year={1994},
  publisher={Wiley Online Library}
}

@article{willems2005note,
  title={A note on persistency of excitation},
  author={Willems, Jan C and Rapisarda, Paolo and Markovsky, Ivan and De Moor, Bart LM},
  journal={Systems \& Control Letters},
  volume={54},
  number={4},
  pages={325--329},
  year={2005},
  publisher={Elsevier}
}

@book{bhatia2013matrix,
  title={Matrix Analysis},
  author={Bhatia, Rajendra},
  year={2013},
  publisher={Springer Science \& Business Media}
}

\newpage
	 \begin{wrapfigure}{l}{1.8cm}
	\includegraphics[width=1.1in,height=1.25in,clip,keepaspectratio]{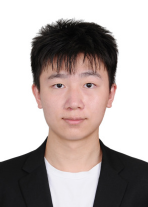}
		 \end{wrapfigure} 
	 Zhaohua Yang received the B.Eng. degree from Zhejiang University, Hangzhou, China, and the B.S. degree from University of Illinois at Urbana-Champaign, Urbana, USA, in 2022. He is currently pursuing the Ph.D. degree with the Department of Electronic and Computer Engineering at The Hong Kong University of Science and Technology, Hong Kong, China. From August 2021 to January 2022, he was a visiting student at the Grainger College of Engineering, University of Illinois at Urbana-Champaign, Urbana, USA.
	 His research interest includes optimal control, data-driven control, and sparse optimization.
	
	 \begin{wrapfigure}{l}{1.8cm}
		 \includegraphics[width=1.1in,height=1.25in,clip,keepaspectratio]{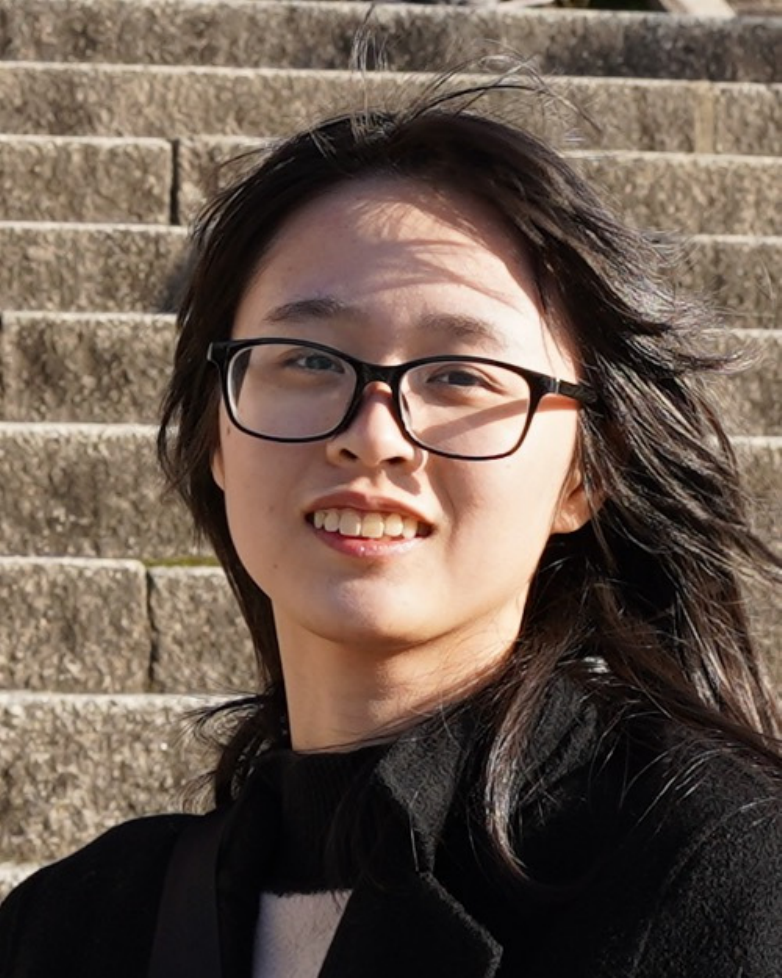}
		 \end{wrapfigure} 
	 Yuxing Zhong received her B.E. degree in Electronic Information and Communications from Huazhong University of Science and Technology, Wuhan, Hubei, China, in 2020, and her Ph.D. degree in Electrical and Computer Engineering from Hong Kong University of Science and Technology (HKUST), Kowloon, Hong Kong in 2024. She is currently a postdoctoral fellow at HKUST. Her research interests include networked control systems, sparsity and event-based state estimation.
	 
	 \begin{wrapfigure}{l}{1.8cm}
		 \includegraphics[width=1.1in,height=1.25in,clip,keepaspectratio]{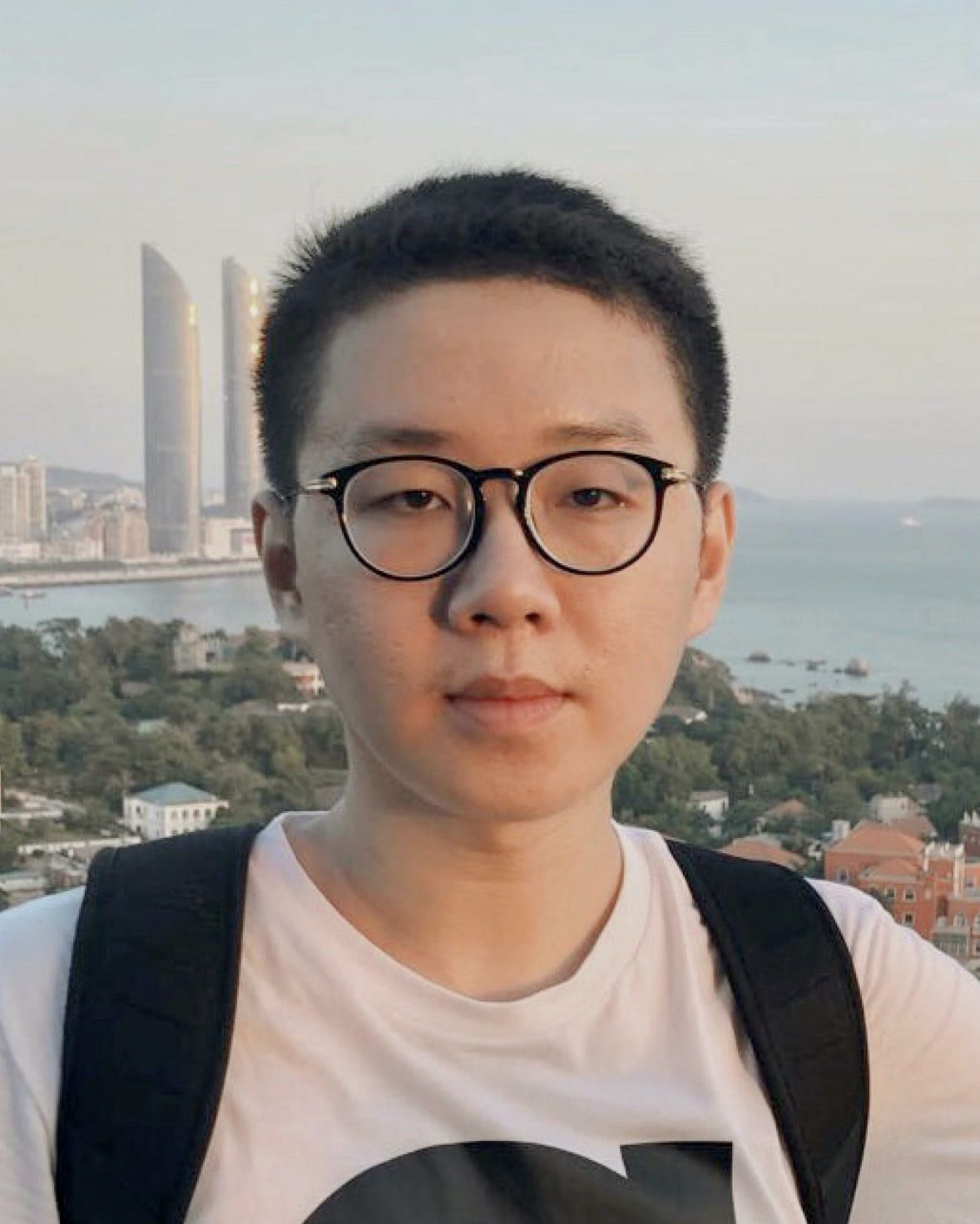}
		 \end{wrapfigure} 
	 Nachuan Yang (Member, IEEE) received his B.Sc. in Mathematics from the University of Hong Kong in 2020 and his Ph.D. in Electronic and Computer Engineering from the Hong Kong University of Science and Technology (HKUST) in 2024. He is currently a Professor with the State Key Laboratory of Environment Characteristics and Effects for Near Space, Beijing Institute of Technology. He serves as an Associate Editor for IEEE Canadian Journal of Electrical and Computer Engineering. His research focuses on learning, optimization, and control of cyber-physical systems.
	
	 \begin{wrapfigure}{l}{1.8cm}
		 \includegraphics[width=1.1in,height=1.25in,clip,keepaspectratio]{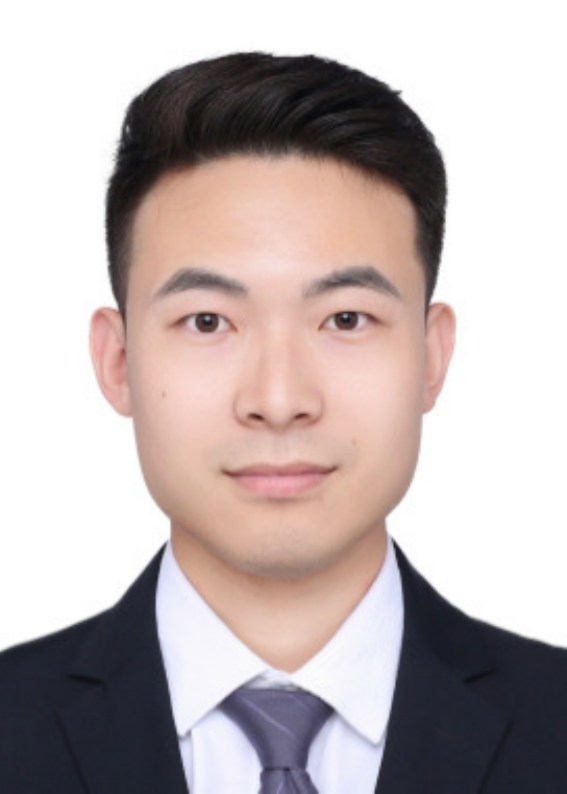}
		 \end{wrapfigure} 
	 Xiaoxu Lyu received his B.Eng. degree in Naval Architecture and Marine Engineering from Harbin Institute of Technology in 2018 and his Ph.D. degree in Dynamical Systems and Control from Peking University in 2023. He is currently a Research Assistant Professor in the Department of Electronic and Computer Engineering at The Hong Kong University of Science and Technology (HKUST). From 2024 to 2025, he was a Postdoctoral Fellow at HKUST. His research interests include distributed estimation and control, data-driven systems, and multi-robot systems.
\newpage
	 \begin{wrapfigure}{l}{1.8cm}
		 \includegraphics[width=1.1in,height=1.25in,clip,keepaspectratio]{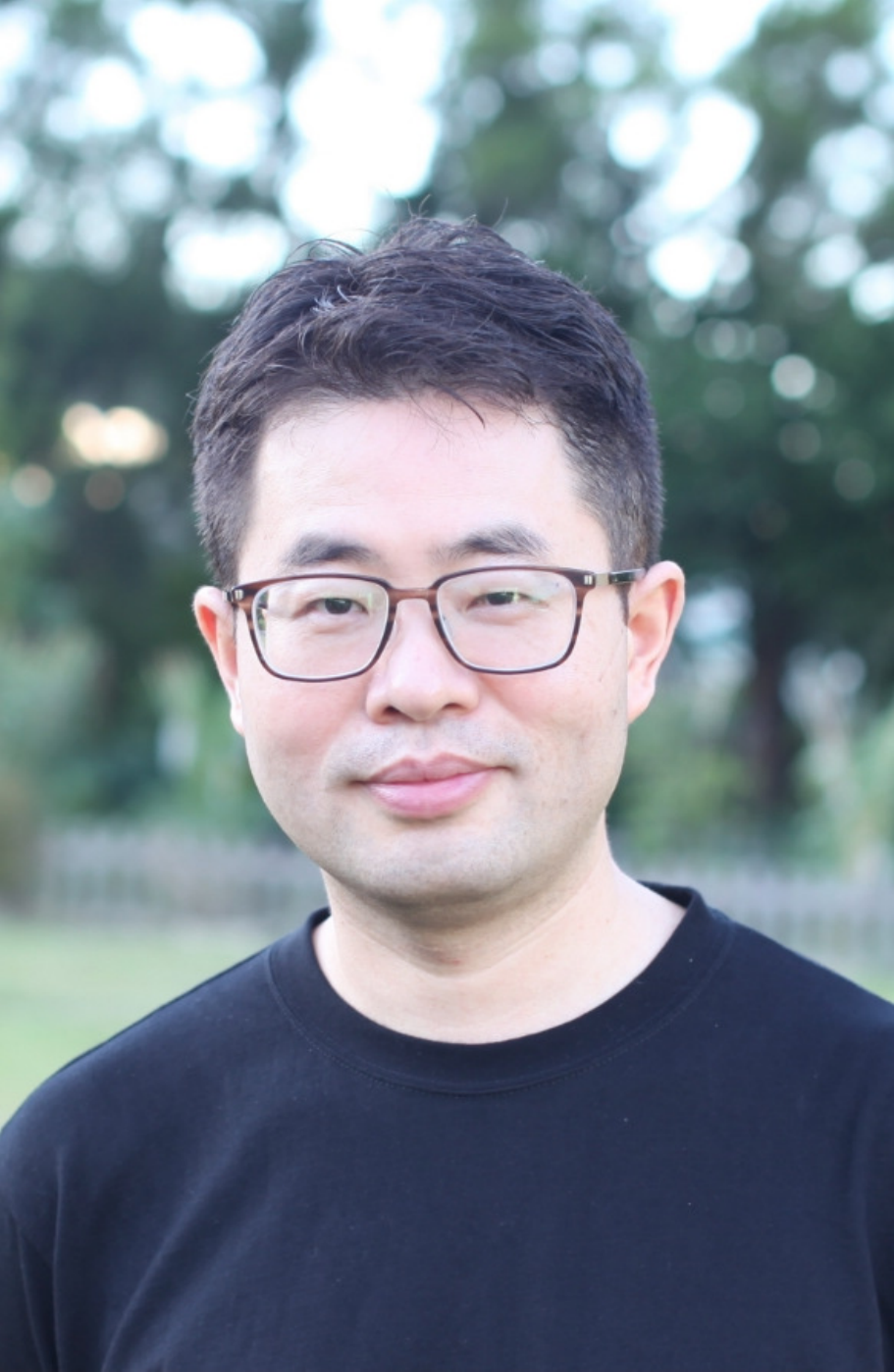}
		 \end{wrapfigure} 
	 Ling Shi received his B.E. degree in Electrical and Electronic Engineering from The Hong Kong University of Science and Technology (HKUST) in 2002 and  Ph.D. degree in Control and Dynamical Systems from The California Institute of Technology (Caltech) in 2008. He is currently a Professor in the Department of Electronic and Computer Engineering at HKUST with a joint appointment in the Department of Chemical and Biological Engineering (2025-2028), and the Director of The Cheng Kar-Shun Robotics Institute (CKSRI). His research interests include cyber-physical systems security, networked control systems, sensor scheduling, event-based state estimation, and multi-agent robotic systems (UAVs and UGVs). He served as an editorial board member for the European Control Conference 2013-2016. He was a subject editor for International Journal of Robust and Nonlinear Control (2015-2017), an associate editor for IEEE Transactions on Control of Network Systems (2016-2020), an associate editor for IEEE Control Systems Letters (2017-2020), and an associate editor for a special issue on Secure Control of Cyber Physical Systems in IEEE Transactions on Control of Network Systems (2015-2017). He also served as the General Chair of the 23rd International Symposium on Mathematical Theory of Networks and Systems (MTNS 2018). He is currently serving as a member of the Engineering Panel (Joint Research Schemes) of the Hong Kong Research Grants Council (RGC) (2023-2026). He received the 2024 Chen Han-Fu Award given by the Technical Committee on Control Theory, Chinese Association of Automation (TCCT, CAA). He is a member of the Young Scientists Class 2020 of the World Economic Forum (WEF), a member of The Hong Kong Young Academy of Sciences (YASHK), and he is an IEEE Fellow.
\end{document}